\documentclass[12pt]{amsart}
\renewcommand{\bar}{\overline}
\usepackage{amsmath,amsthm,amscd,euscript}
\def \r{\mathbb R}

\def \z{\mathbb Z}

\DeclareMathOperator{\isin}{lsin}

\DeclareMathOperator{\il}{l\ell}

\newtheorem{theorem}{Theorem}

\newtheorem{proposition}[theorem]{Proposition}
\newtheorem{corollary}[theorem]{Corollary}
\theoremstyle{remark}

\theoremstyle{definition}
\newtheorem{definition}[theorem]{Definition}
\newtheorem{example}{Example}
\newtheorem{problem}{Problem}

\title{On determination of periods of geometric continued fractions
for two-dimensional algebraic hyperbolic operators.}
\email[Oleg Karpenkov]{karpenk@mccme.ru}

\begin{document}
\input epsf


\sloppy \normalsize


\begin{center}
{\bf On determination of periods of geometric continued fractions
for two-dimensional algebraic hyperbolic operators.} \vspace{5mm}

{\sc O.~N.~Karpenkov\footnote[1]{The work is partially supported
by RFBR SS-1972.2003.1 grant, by RFBR grant 05-01-01012a, by
NWO-RFBR 047.011.2004.026 (RFBR 05-02-89000-NWO\_a) grant,
NWO-DIAMANT 613.009.001 grant, and by RFBR 05-01-02805-CNRSL\_a
grant.} }

Mathematisch Instituut, Universiteit Leiden, P.O. Box 9512, 2300
RA Leiden, The Netherlands)
\end{center}


\section*{Introduction}

A two-dimensional operator $A$ in $SL(2,\z)$ is said to be {\it
hyperbolic}, if its eigenvalues are real and distinct. In the
present paper we study the connection between periods of
geometric continued fraction in the sense of Klein and reduced
operators, described by J.~Lewis and D.~Zagier in~\cite{DZg}.
Actually, a determination of a period of the geometric continued
fraction for the operator $A$ is equivalent to a calculation of a
period of the ordinary continued fraction for the tangent of an
angle between any eigen straight line of $A$ and the horizontal
axis.

In present paper for any period of a geometric continued fraction
in the sense of Klein we make an explicit construction of a
reduced hyperbolic operator in $SL(2,\z)$ with the given period
for the geometric continued fraction (Theorem~\ref{th}). It turns
out that reduced operators naturally forms one-parametric
families. Further we describe an algorithm to construct a period
for an arbitrary operator of $SL(2,\z)$. The base part of the
algorithm is to determine a reduced operator that is conjugate to
the given (the Gauss Reduction Theory).

In 1993 V.~I.~Arnold formulated some questions on periods of
continued fractions related to the eigenvectors of the
$SL(2,\z)$-operators, see for instance in~\cite{ArnProb}
and~\cite{Arn2}. The first studies of these problems are made in
the article~\cite{MPav} by M.~Pavlovskaya, in which the author
experimentally investigates statistical questions on geometrical
continued fractions (such as verification of the average length of
periods of continued fractions). The questions on the
distribution of positive integers in minimal periods of quadratic
continued fractions were studied by M.~O.~Avdeeva and
B.~A.~Bykovskii in the works~\cite{Avd1} and~\cite{Avd2}. In his
work~\cite{ArnAr} V.~I.~Arnold investigates palindromic
properties of continued fraction periods and the connection to
the integer forms (the complete proof of the palindromic property
is given by F.~Aicardi and M.~Pavlovskaya). The relation between
``$-$''-continued fractions of hyperbolic operators and Fuchian
groups, and a few words on the algorithm of integer conjugacy
check for the operators in $SL(2,\z)$ is described in the work of
S.~Katok~\cite{Kat}. Some estimates of the period lengths for
ordinary continued fractions for $\sqrt{d}$ are given
in~\cite{Hic} by D.~R.~Hickerson.

In the works~\cite{Kle1} and~\cite{Kle2}  F.~Klein generalized
geometric continued fraction to the multidimensional case of
$SL(n,\z)$-operators. For more information see the works of
V.~I.~Arnold~\cite{Arn2}, E.~I.~Korkina~\cite{Kor2},
G.~Lachaud~\cite{Lac1}, J.-O.~Moussafir~\cite{Mou2}, the
author~\cite{Kar1}, etc.

The work is organized as follows. In the first section we give
definitions of ordinary continued fractions and geometric
continued fractions in the sense of Klein, we show a duality of
sails of a geometric continued fraction. Further in the second
section we construct families of operators possessing the
geometric continued fractions (actually the corresponding
LLS-sequences) with the period $(a_1,\ldots,a_{2n-1},t)$, where
$t$ is the parameter of families. Further we give an algorithm of
period calculation for geometric continued fractions. In the
third section we formulate some questions and show the results of
a few experiments related to the algorithm.

The author is grateful to V.~I.~Arnold for constant attention to
this work
and Mathematisch Instituut of Universiteit Leiden for the
hospitality and excellent working conditions.

\section{Geometric continued fractions in the sense of Klein}

{\bf Ordinary continued fractions.} Consider an arbitrary finite
sequence of integers $(a_0, a_1, \ldots ,a_n)$, where $a_0$ is
any integer and $a_i>0$ for $i>0$. The rational
$$
a_0+1/(a_1+1/(a_2+\ldots)\ldots))
$$
is called the {\it ordinary continued fraction} and denoted by
$[a_0: a_1; \ldots ;a_n]$.

A continued fraction with even (odd) number of elements is called
an {\it even $($odd$)$} ordinary continued fraction respectively.

\begin{proposition}
For any rational there exist a unique even and a unique odd
continued fractions. \qed
\end{proposition}

For instance the even and odd ordinary continued fractions for a
rational $47/39$ are $[1{:}4;1;7]$ and $[1{:}4;1;6;1]$
respectively.

\vspace{3mm}

{\bf Sails of hyperbolic operator.} By $[[a,c][b,d]]$ we denote
the operator
$$\left(
\begin{array}{cc}
a&c\\
b&d\\
\end{array}
\right).
$$

An operator in $SL(2,\r)$ with two distinct real eigenvalues is
called {\it hyperbolic}.

A vector (segment) is said to be {\it integer} if it (its
endpoints) has integer coordinates.

Consider a hyperbolic operator $A$ with no integer eigenvectors.
The operator $A$ has exactly two distinct eigen straight lines.
These lines do not contain integer points of the lattice distinct
to the origin. The complement to the union of the lines consists
of four piece-wise connected components, each of which is an open
octant. Consider one of these octants. The boundary of the convex
hull of all integer points except the origin in the closure of
the octant is called a {\it sail} of the operator $A$. The set of
all sails is called  the {\it geometric continued fraction in the
sense of Klein} for the operator $A$ (see also the works of
F.~Klein~\cite{Kle1}, E.~I.~Korkina~\cite{Kor2}, and
V.~I.~Arnold~\cite{Arn2}). Two sails are said to be {\it
equivalent} if there exists a linear lattice-preserving
transformation, taking one of the sails to another. Geometric
continued fractions with equivalent sails are called {\it
equivalent}.

\vspace{1mm}

{\it Remark.} The majority of constructions of this article can
be naturally generalized to the case of operators that have
integer eigenvectors. For simplicity of the exposition we do not
consider such operators in the article.

\vspace{3mm}

{\bf LLS-sequence.} An {\it integer length} of an integer segment
$PQ$ is the quantity of inner integer points of the segment plus
one, it is denoted by $\il(PQ)$. Let integer segments $PQ$ and
$PR$ do not lie on the same straight line. An {\it integer sine}
of the angle $QPR$ is an index of the sublattice generated by the
integer vectors of the straight lines $PQ$ and $PR$ in the
lattice of all integer vectors (we denote it by $\isin (QPR)$).
An integer sine of an angle contained in some straight line is
supposed to be equivalent to zero.

Any sail of the hyperbolic operator $A$ with no integer eigen
vectors is a two-sides broken line with infinite in both sides
number of segments.

\begin{definition}
Let $\ldots V_{-2}V_{-1}V_0V_1V_2 \ldots$ be a sail of some
operator. The infinite in both sides sequence of positive integers
$$
(\ldots,\il(V_{-2}V_{-1}),\isin\angle V_{-2}V_{-1}V_0,
\il(V_{-1}V_{0}),\isin\angle V_{-1}V_{0}V_1,
\il(V_{0}V_{1}),\isin\angle V_{0}V_{1}V_2, \ldots )
$$
is called the {\it LLS-sequence} of the sail (see also
in~\cite{KarPrep}).
\end{definition}

Two LLS-sequences are called {\it equivalent}, if one can be
obtained from another by shifting on a finite number of elements
and/or reversing the order of all elements. The LLS-sequences of
equivalent sails are equivalent, since all integer lengths and
sines are invariants under integer-linear transformations.

\begin{figure}
$$\epsfbox{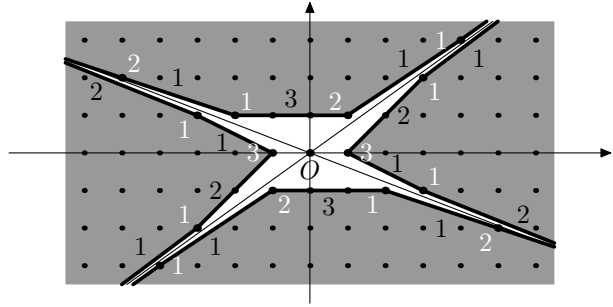}$$
\caption{Geometric continued fraction of the operator
$[[7,18][5,13]]$.}\label{cf.1}
\end{figure}

On Figure~\ref{cf.1} we show the geometric continued fraction for
the operator $[[7,18][5,13]]$. Integer lengths of edges are
denoted by black digits, and integer sines --- by white. The
LLS-sequences of all four sails are equivalent to
$(\ldots,2,1,1,3,2,1,1,\ldots)$.

\vspace{3mm}

{\bf Duality of sails.} Two sails are called {\it dual} with
respect to each other if their LLS-sequences coincide, the
sequence of integer lengths of the first sail coincides with the
sequence of integer sines of the second, and the sequence of
integer sines of the first sail coincides with the sequence of
integer lengths of the second. (On the duality for
multidimensional continued fractions in the sense of Klein see in
the book by G.~Lachaud~\cite{Lac1}.)

\begin{proposition}\label{st1}
Let $A$ be a hyperbolic operator with no integer eigenvectors.
Then the LLS-sequences of all his four octants coincide. Moreover
the sails of the opposite octants are equivalent. The sails of the
adjacent octants are dual.
\end{proposition}

In the proof of Proposition~\ref{st1} we use the following
properties of an integer sine (for further information see
in~\cite{KarPrep}):
$$
\isin(PQR)=\isin(RQP); \quad \isin(PQR) = \isin (\pi -PQR),
$$
where by $\pi {-}PQR$ we denote the angle adjacent to the angle
$PQR$.

The statement of Proposition~\ref{st1} is known to the author from
V.~I.~Arnold and E.~I.~Korkina and supposed to be classical, still
its formulation and proof seem to be missing in the literature.
To avoid further questions we give the proof here.

\begin{proof}
The sails of the opposite octants (and so their LLS-sequences)
are equivalent, since they are symmetric about the origin.

Consider now the case of adjacent sails. Let now $S$ and $S'$ be
adjacent sails of the operator $A$, where $S'$ is counterclockwise
with respect to $S$ while turning around the origin $O$. Let $S$
be a broken line $\ldots V_{-1}V_0V_{1}\ldots$ with
counterclockwise order of vertices with respect to the origin.

Let us show a natural bijection between the edges of the sail $S$
and vertices of the sail $S'$. Consider an edge $V_{i}V_{i+1}$ of
the sail $S$. Denote by $V$ the closest integer point in the
segment $V_{i}V_{i+1}$  to the endpoint $V_{i+1}$ and distinct to
$V_{i+1}$. Denote the integer point $O{+}\bar{VV_{i+1}}$ by
$V_i'$. Notice that the point $V_i'$ is in the convex hull of
$S'$.

Since the triangle $OVV_{i+1}$ does not contain integer points
inside, the triangle $OV_i'V_{i+1}$ does not contain integer
points inside. Hence there is no integer points between the
parallel lines $OV$ and $V_i'V_{i+1}$. Therefore, $V_i'$ is a
vertex of the sail $S'$ (see on Figure~\ref{proof.123}a).

\begin{figure}
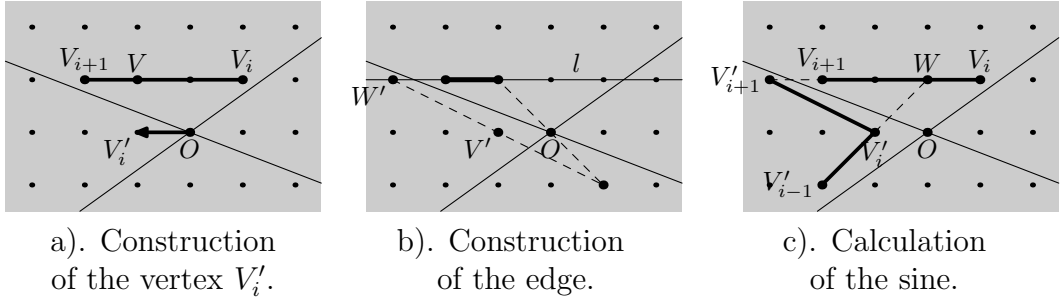

$$
\begin{array}{ccc}
\epsfbox{proof.1}&\epsfbox{proof.2}&\epsfbox{proof.3}\\
\hbox{a). Construction}&\hbox{b). Construction }&\hbox{c). Calculation}\\
\hbox{of the vertex $V'_i$.}&\hbox{of the edge.}&\hbox{of the sine.}\\
\end{array}
$$
\caption{Duality of adjacent angles.}\label{proof.123}
\end{figure}

From the other side $V'$ is a vertex of the sail $S'$. Denote by
$l$ the closest to the line $OV'$ parallel line containing
integer points and intersecting $W$. The line $l$ intersects with
the octant containing $S'$ in a ray. Denote by $W'$ the closest to
the vertex of the ray integer point of this ray. Then the point
$V'{+}\bar{W'V'}$ is in the octant opposite to the octant with
the sail $S$ (otherwise $V'$ is not a vertex). Hence the point
$W'{+}2\bar{V'O}$ symmetric about the origin to the point
$V'{+}\bar{W'V'}$ lies in the octant of the sail $S$. Therefore,
the integer points $W'{+}\bar{V'O}$ and $W'{+}2\bar{V'O}$ lie in
the octant containing $S$. Since $l$ is the closest to the line
$XV'$ parallel line containing integer points and intersecting
$S$, the segment with the endpoints $W'{+}\bar{V'O}$ and
$W'{+}2\bar{V'O}$ in contained in the sail (see on
Figure~\ref{proof.123}b).

Therefore the above correspondence between the edges
$V_{i}V_{h+2}$ of the sail $S$ and vertices $V_i'$ of the sail
$S'$ is a bijection. Moreover the order of vertices $V_i'$ is
clockwise.

Let us prove now that
$\il(V_iV_{i+1})=\isin(V_{i-1}'V_i'V_{i+1}')$. Denote the point
$V_i{+}\bar{OV_i'}$ by~$W$. By construction, the point $W$
belongs to the segment $V_iV_{i+1}$ (see
Figure~\ref{proof.123}c). Then
$$
\isin(V_{i-1}'V_i'V_{i+1}')=\isin(\pi-WV_i'V_{i+1}')=\isin(WV_i'V_{i+1}')
=\il(WV_{i+1}')=\il(V_iV_{i+1}).
$$

By the analogous argumentation the following equalities also hold:
$$
\il(V_i'V_{i+1}')=\isin(V_{i-1}V_iV_{i+1}).
$$

Therefore, the LLS-sequences of the sails $S$ and $S'$ coincide,
moreover the sequence of integer lengths (sines) of $S$ coincides
with the sequence of integer sines (lengths) of $S'$. Hence $S$
and $S'$ are dual. Proposition~\ref{st1} is proved.
\end{proof}

\vspace{3mm}

{\bf Existence and uniqueness of the equivalence classes of
continued fractions with a given LLS-sequences.}

\begin{definition}
The {\it LLS-sequence} of an operator $A$ is the LLS-sequence for
any of its sails.
\end{definition}

\begin{proposition}~\label{st2}
i$)$. For any infinite in two sides sequence of integers there
exists a hyperbolic operator whose LLS-sequence coincides with the
given.

ii$)$. Two sails with coinciding LLS-sequences are either
equivalent or dual. All sails dual to the given are equivalent to
each other.
\end{proposition}

\begin{proof}
Both statements follows directly from Corollary~5.11
of~\cite{KarPrep} (see also in the work of
E.~I.~Korkina~\cite{Kor2}).
\end{proof}

\section{Algebraic sails.}

{\bf Construction of the operator with the given period.}
Consider now an algebraic case of hyperbolic operators of the
group $SL(2,\z)$ with irreducible over rationals characteristic
polynomial. Let $A$ be an integer hyperbolic operator of
$SL(2,\z)$. Denote by $\Xi(A)$ the group of operators in all
$SL(2,\z)$ commuting with $A$ and having positive real
eigenvalues. By Dirichlet unity theorem~\cite{BSh} the group
$\Xi(A)$ is isomorphic to $\z$. Any sail of the operator $A$ is
invariant under the action of the group $\Xi(A)$, moreover the
operators of the group $\Xi(A)$ act on the sails by shifting the
edges of the broken line along thy broken line. The sails of the
operator $A$ are called {\it algebraic}.

Therefore LLS-sequences of hyperbolic algebraic operators are
periodic. The converse is also true (see Corollary~\ref{st3}
below).

\vspace{1mm}

{\it Remark.} On Figure~\ref{cf.1} we show the sails of
hyperbolic algebraic operator $[[7,18][5,13]]$ with a period of
the LLS-sequence equals $(2,1,1,3)$.

\vspace{1mm}

\begin{theorem}\label{th}
Consider an $ST(2,\z)$-operator $[[a,c+\lambda a][b,d+\lambda
b]]$.
\\
i$)$. Let $a=0$, $b=1$, $d=1$. If $\lambda>2$ then the operator
$A$ is hyperbolic and its sails are algebraic. One of the periods
of the LLS-sequence of the operator  $A$ is
$$
(1,\lambda-1).
$$
ii$)$. Let $b>a \ge 1$, $0<d\le b$, $\lambda\ge 1$. Let the odd
ordinary continued fraction for $b/a$ equal
$$
[a_0{:}a_1;\ldots;a_{2n}].
$$
Then the operator $A$ is hyperbolic and its sails are algebraic.
One of the periods of the LLS-sequence of the operator  $A$ is
$$
(a_0,a_1,\ldots,a_{2n},\lambda).
$$
\end{theorem}

Notice that for any couple of relatively prime integers $(a,b)$
where $b>a \ge 0$ there exists a couple of integers $(c,d)$,
satisfying $0<d\le b$ and $ad{-}bc=1$.

\vspace{1mm}

{\it Remark.} For negative values of $\lambda$ in the case $a=0$
the periods are $(1,|\lambda|{-}3)$. In the case $a>0$ the
periods equal
$$
(a'_0,a'_1,\ldots,a'_{2m},|\lambda|-2),
$$
where $[a'_0{:}a'_1;\ldots;a'_{2m}]$ --- is the odd ordinary
continued fraction for $b/(b-a)$.

\vspace{1mm}

\begin{proof}
The discriminant of the characteristic polynomial of the operator
$A$ equals $((a+\lambda b+d)^2-4$. Since $\lambda\ge 1$, $b\ge
1$, $d\ge 1$, and $a\ge 0$, the discriminant is nonnegative.
Besides it equals zero in the exceptional case $a=0$, $b=1$,
$\lambda=1$, $d=1$. Therefore the operator $A$ is hyperbolic in
all cases. Since $t^2-4$ for integer $t>2$ is not a square of some
integer, the sails of the operator $A$ are algebraic.

Let us now construct a period for the LLS-sequence. Note that
both eigenvalues of the operator $A$:
$$
\frac{a+\lambda b+d\pm \sqrt{((a+\lambda b+d)^2-4}}{2}
$$
are positive, and thus the operator $A$ takes each sail to
itself. Consider the sail $S$, whose convex hull contains the
point $P=(1,0)$.

\vspace{1mm}

Suppose the operator $A$ is of series i). Then the set of the
vertices for one of its sails coincides with the set of points
$A^n(1,0)$ with an integer parameter $n$. Simple calculations
lead to the result of the theorem.

\vspace{1mm}

Let now an operator $A$ be an operator of series ii), i.~e.
$b>a>0$.

Denote by $\alpha$ the closed convex angle with vertex at the
origin and edges passing through the points $P$ and $A(P)$.
Consider the boundary of the convex hull of all integer points
inside $\alpha$ except the origin. The boundary consists of two
rays and a finite broken line. We call the finite broken line in
the boundary the {\it sail of the angle $\alpha$} and denote it by
$S_\alpha$.

Now we show that the sail of the angle $\alpha$ is completely
contained in one of the sails of the operator $A$. Denote by
$S_\alpha^\infty$ the following infinite broken line:
$$
\bigcup\limits_{i=-\infty}^{+\infty} \Big(A^i(S_\alpha)\Big).
$$
By the construction the convex hull of $S_\alpha^\infty$ coincides
with the convex hull of the sail $S$. It remains to verify if
$S_\alpha^\infty$ coincides with the boundary of its convex hull.
The broken line $S_\alpha^\infty$ is the boundary of the convex
hull if the convex angles generated by adjacent edges of the
broken line do not contain the origin. To check this it is
sufficient to study all the angles of one of the periods of the
broken line $S_\alpha^\infty$, for instance all the angles with
vertices at vertices of $S_\alpha$ except the point $A(P)$. All
convex angles generated by couples of adjacent edges at inner
vertices of the sail $S_\alpha$ do not contain the origin by
definition. It remains to check the angle with vertex at
$P=(1,0)$.

Since $b/a>1$, the first edge is parallel to the vector $(0,1)$.
Consider the second edge $PQ$. Since the triangle $OPQ$ does not
contain integer points distinct to the integer points of the
segment $PQ$ and the vertex $O$, the segment $PQ$ contains the
point with coordinates $(x,-1)$. By the convexity of finite
broken line $A^{-1}(S_\alpha)$ the value of $x$ is determined by
the eigen direction:
$$
\left(\frac{-a+\lambda b+d+ \sqrt{((a+\lambda
b+d)^2-4}}{2b},-1\right),
$$
namely,
$$
x=\left\lfloor\frac{-a+\lambda b+d+ \sqrt{((a+\lambda
b+d)^2-4}}{2b}\right\rfloor+1,
$$
where $\lfloor t\rfloor$ denotes the maximal integer that does not
exceeding $t$. As one can show, $x$ is contained in the open
interval $\big( \lambda{+}1{+}(d{-}1)/b,\lambda{+}1{+}d/b \big)$.
By condition $b\ge d>0$ we have
$$
x=\lambda+1.
$$
Hence for $\lambda>1$ the convex angle at vertex $P$ does not
contain the origin. Therefore, the broken line $S_\alpha^\infty$
coincides with the sail.

In the paper~\cite{KarPrep} it is shown that the sail $S_\alpha$
consists of $n{+}1$ segments. The integer lengths of the
consecutive segments equal $a_0,a_2,\ldots, a_{2n}$, and the
integer sines of the corresponding angles equal $a_1,a_3,\ldots,
a_{2n-1}$ respectively. Now note, that from the explicit value of
$x$ it follows that the integer sine for the angle at point $P$
equals $\lambda$. Hence the LLS-sequence of the sail $S$ has a
period
$$
(a_0,a_1,\ldots,a_{2n},\lambda).
$$
Therefore the LLS-sequence of the operator $A$ has the prescribed
period.
\end{proof}

\begin{corollary}\label{st3}
A sail with the periodic LLS-sequence is algebraic $($i.~e. a
sail of some algebraic hyperbolic operator$)$.
\end{corollary}

\begin{proof}
In Theorem~\ref{th} we constructed the algebraic operators for
all finite sequences as periods. Then in Proposition~\ref{st2} we
showed that the sails with equivalent LLS-sequences are either
equivalent or dual. Therefore any sail with periodic LLS-sequence
is algebraic.
\end{proof}

\vspace{1mm}

{\it Remark.} Consider some sail with periodic LLS-sequence. Let a
minimal period of LLS-sequence is even and consists of $2n$
elements. Then there exists an $SL(2,\z)$-operator $A$  with
positive eigenvalues, that makes an $n$-edge shift of the sail
along the sail. Precisely this operator generates the group
$\Xi(A)$ of the sail shifts (see above). Let a minimal period of
LLS-sequence is odd and consists of $2n{+}1$ elements (in
particular this implies that the sail is equivalent to any dual
sail). Then there exists a $GL(2,\z)$-operator $B$ with negative
discriminant, whose square makes an $(2n{+}1)$-edge shift of the
sail along the sail. Moreover, the operator $B^2$ generates the
group $\Xi(T)$.

\vspace{1mm}

{\it Remark.} Let us say a few words about non-hyperbolic
operators in $SL(2,\z)$. It turns out that each of such operators
is equivalent to exactly one of the operators of the following
list:

---  $[[1,1][-1,0]]$;

---  $[[0,1][-1,0]]$;

---  $[[0,1][-1,-1]]$;

---  $[[1,n][0,1]]$, where $n$ is any integer.

\vspace{3mm}

{\bf Algorithm of finding a period of the LLS-sequence for a
hyperbolic algebraic operator.}

\begin{definition}
An operator $[[a,c][b,d]]$ in $SL(2,\z)$ is said to be {\it
reduced}, if the following holds: $d> b\ge a\ge 0$.
\end{definition}

{\it Remark.} The definition of a reduced operator is slightly
different to one given in the works~\cite{DZg} and~\cite{Man}:
{\it an operator in $SL(2,\z)$ is reduced iff it has non-negative
entries which are non-decreasing downwards and to the right}.

The main idea of the calculation of the period is to find a
reduced operator with the sails equivalent to the sails of the
given one. Then it remains to calculate the period of the reduced
operator by Theorem~\ref{th}.

\vspace{1mm}

{\it Data.} Suppose we know the integer entries of an operator
$[[a,c][b,d]]$ with unit determinant and positive discriminant.
Let also the characteristic polynomial does not has roots $\pm 1$.
From the listed conditions it follows that the operator
$[[a,c][b,d]]$ is hyperbolic operator in $SL(2,\z)$ with
irreducible characteristic polynomial.

\vspace{1mm}

{\it It is requested} to construct one of the periods of the
LLS-sequence of the hyperbolic algebraic operator $[[a,c][b,d]]$.

\vspace{1mm}

{\bf Description of the algorithm.}

{\it Step 1.} If $b<0$, then we multiply the operator
$[[a,c][b,d]]$ by $[[-1,0][0,-1]]$. The LLS-sequence does not
change at that.

{\it Step 2.} So, now the element $b$ is positive. Conjugate the
operator $[[a,c][b,d]]$ by the operator $[[1,-\lfloor a/b
\rfloor][0,1]]$. We obtain the operator $[[a',b'][c',d']]$, where
$0\le a' \le b'$.

{\it Step 3.1.} Suppose $b'=1$, then $a'=0$, $c'=-1$. Moreover we
have $|d|>2$, since otherwise the original operator is not
algebraic. Therefore a period of the LLS-sequence equals
$(3,|d|-2)$.

{\it Step 3.2.1.} Suppose, $b'>1$. If $d'>b'$, then we have found
a reduced operator, now we go to Step~4.

{\it Step 3.2.2.} Suppose, $b'>1$. If $d'<-b'$, then we conjugate
by the operator $[[-1,1][0,1]]$. Finally we have the operator
$[[a'',c''][b'',d'']]$ with $b''=b'$, $a''=b'-a'$, and
$d''=-b'-d'>0$, further we should go to Step~3.2.1, or to
Step~3.2.3.

{\it Step 3.2.3.} Suppose, $b'>1$. The case $|d'|\le |b'|$. Note
that the absolute values of $b'$ and $d'$ do not coincide since
the determinant of the operator does not have divisors distinct
to the unity. Therefore it remains the case $|d|<|b|$. In this
case we have:
$$
|c'|=\left|\frac{a'd'-1}{b'}\right|\le\frac{(b'-1)^2+1}{b'}\le
b'-1.
$$
We conjugate the operator $[[a',c'][b',d']]$ with the operator
$[[0,-1][-1,0]]$ and obtain $[[d',b'][c',a']]$, where $|c'|<|b'|$.
Now we return back to Step~1 with the obtained operator
$[[d',b'][c',a']]$.

{\it Step 4.} We obtained a reduced operator $[[\tilde a,\tilde
c][\tilde b,\tilde d]]$, $\tilde b > 1$ with the LLS-sequence
equivalent to the LLS-sequence of the original operator. By
Theorem~\ref{th} to construct one of the periods of the
LLS-sequences of the reduced operator we should construct the odd
ordinary continued fraction for $\tilde b/\tilde a$, and find the
integer $\big\lfloor(\tilde d{-}1)/\tilde b \big\rfloor$.

\section{Some questions and examples}

{\bf A question on complexity of the minimal period.} Note that
for any operator there exist finitely many reduced operators with
the same trace and LLS-sequence. If we study the reduced operators
that make shifts of sails on a minimal possible period, then the
number of such operators coincides with the length of the minimal
period (see also in~\cite{DZg}). Let $[[a,c][b,d]]$ be a reduced
operator. We call the integer $b$ --- its {\it complexity}.

\begin{problem}
Study the minimal complexity for reduced operators with
LLS-sequence having a length $n$ period $(a_1,\ldots, a_n)$.
\end{problem}

{\it Remark.} The minimal complexity coincides with the minimal
positive value of the integer sine of the angles $POQ$, where $O$
is the origin, $P=(x,y)$ is an arbitrary integer point distinct to
$O$, and $Q=A(P)$. Therefore the minimal complexity, considered
as the minimal possible integer sine, is well defined for all
operators and it is invariant under conjugations.

If $n$ is even, then the problem is equivalent to finding the
minimal numerator among the numerators of the rationals:
$$
[a_1{:}\ldots;a_{n-1}], \quad [a_2{:}\ldots;a_{n}], \quad
[a_3{:}\ldots;a_{n};a_1], \quad \ldots \quad
,[a_n{:}a_1;\ldots;a_{n-2}].
$$

\begin{example} Let the period contains two elements: $(a,b)$, where $a<b$,
then the minimal complexity equals $a$.
\end{example}

\begin{example}
Let the period contains four elements: $(a,b,c,d)$. Let $d$ is
not smaller than the other elements of the period, let also $d>a$
except for the case $a=b=c=d$. Then the rational with the minimal
numerator can be found from the following table.

\begin{center}
\begin{tabular}{|c||c|}
\hline
$(a,b,c,d)$    & Rational  with          \\
$d\ge a,b,c$  &the minimal numerator     \\
\hline \hline
$d> a,b,c$       &$[a{:}b;c]$             \\
$d=c;b<a<d$      &$[d{:}a;b]$             \\
$d=c;a<b<d$      &$[a{:}b;c]$             \\
$d=c;a=b<d$      &$[a{:}b;c]$ and $[d{:}a;b]$ \\
$d=b;d>a,c$      &$[a{:}b;c]$ and $[c{:}d;a]$ \\
$d=c=b;d>a$      &$[a{:}b;c]$ and $[c{:}d;a]$ \\
$d=c=b=a$        &all                   \\
\hline
\end{tabular}
\end{center}
\end{example}

If $n$ is odd, then the problem is equivalent to finding the
minimal numerator among the numerators of the rationals (we define
$a_{n+k}=a_k$):
$$
[a_1{:}\ldots;a_{2n-1}], \quad [a_2{:}\ldots;a_{2n}], \quad
[a_3{:}\ldots;a_{2n};a_1], \quad \ldots \quad
,[a_{2n}{:}a_1;\ldots;a_{2n-2}].
$$

\begin{example}
Let the period consists of three elements: $(a,b,c)$, where $c\ge
a,b$. Then the fraction $[a{:}b;c;a;b]$ has the minimal numerator
(or of one of some equivalent minimal numerators in the case of
$a=c$ or $b=c$).
\end{example}

\vspace{1mm}

One can suppose that we should skip one of the maximal elements
of the period, but that is not true for the six element sequence:
$(1,4,5,4,1,4)$. The minimum of the numerators is attained at the
fraction $[1:4;5;4;1]$, and not at the fraction $[4:1;4;1;4]$.

\vspace{3mm}

{\bf On frequencies of occurrences of the reduced operators.}
First we describe a proper probabilistic space. Let
$P=(a_1,a_2,\ldots,a_{2n-1},a_{2n})$ be some period. Denote by
$\#_N(P)$ the quantity of all operators satisfying the following
conditions:\\
{\it i}). The absolute value of any entry of the operator does not
exceed $N$.
\\
{\it ii}). The sequence $P$ is one of the periods of SL-sequence
for the operator.
\\
{\it iii}). Starting from the operator, the algorithm of the
previous section constructs the reduced operator $[[a,b][c,d]]$,
where $(a,b)=(0,1)$ for the case $P=(1,a_2)$; and
$$
b/a=[a_1:a_2;\ldots;a_{2n-1}]
$$
in the remaining cases.

Then one studies the relative statistics of $\#_N(P)$ while $N$
tends to infinity. The following questions are of interest.

\begin{problem}
{\bf a).} Which one of the reduced operators with a given trace
(or with a fixed LLS-sequence) is the most frequent as a result of
the algorithm of the previous section?\\
{\bf b).}  What is the probability of that? \\
{\bf c).} Is it true that the maximal possible probability is
attained at reduced operators with minimal complexity?
\end{problem}

In Table~1 we give some results of calculation of $\#_{25000}(P)$
for the operators with small absolute value of the trace. We
remind that the minimal absolute value of the trace of hyperbolic
$SL(2,\z)$-operator equals $3$.

It is interesting to note that $SL(2,\z)$-operators corresponding
to $P=(1,2)$ are more frequent than the $SL(2,\z)$-operators
corresponding to $P=(1,1)$. This occurs since the sails whose
LLS-sequences has the period $(1,1)$ are equivalent to their
duals. If we enumerate the operators with multiplicities
equivalent to the number of equivalent sails for the operators
then we get:
$$
4\#_{25000}(1,1)>2\#_{25000}(1,2)+2\#_{25000}(2,1).
$$

\begin{table}\label{Tab25000}
\begin{center}
\begin{tabular}{|c|c|c|c|c|}
\hline
Absolute value      &Notation for classes        &Period       & Operator     & Value of       \\
of the trace        &of equivalent operators     &$P$          &$[[a,b][c,d]]$&$\#_{25000}(P)$ \\
\hline \hline
3  & $L_3$& $(1,1)$ & $[[0,1][-1,3]]$ & $663160$\\
\hline
4  & $L_{4}$ &  $(1,2)$ & $[[0,1][-1,4]]$   & $834328$\\
   &         &  $(2,1)$ & $[[1,2][1,3]]$    & $304776$\\
\hline
5  & $L_{5}$ &  $(1,3)$ & $[[0,1][-1,5]]$   & $818200$\\
   &         &  $(3,1)$ & $[[1,3][1,4]]$    & $194528$\\
\hline
6  & $L_{6,1}$& $(1,4)$ & $[[0,1][-1,6]]$   & $777128$ \\
   &          & $(4,1)$ & $[[1,4][1,5]]$    & $141784$ \\
\cline{2-5}
   & $L_{6,2}$& $(2,2)$ & $[[1,2][2,5]]$    & $446432$\\
\hline
7  & $L_{7,1}$& $(1,5)$ & $[[0,1][-1,7]]$   & $734904$\\
   &          & $(5,1)$ & $[[1,5][1,6]]$    & $110848$\\
\cline{2-5}
   & $L_{7,2}$& $(1,1,1,1)$ & $[[2,3][3,5]]$& $201744$\\
\hline
8  & $L_{8,1}$ &$(1,6)$ & $[[0,1][-1,8]]$   & $695560$\\
   &           &$(6,1)$ & $[[1,6][1,7]]$    & $90688$ \\
\cline{2-5}
   & $L_{8,2}$ &$(2,3)$ & $[[1,2][3,7]]$    & $435472$\\
   &           &$(3,2)$ & $[[1,3][2,7]]$    & $310872$\\
\hline
9  & $L_{9}$   &$(1,7)$ & $[[0,1][-1,9]]$   & $660984$\\
   &           &$(7,1)$ & $[[1,7][1,8]]$    & $76552$ \\
\hline
10 & $L_{10,1}$ & $(1,8)$ & $[[0,1][-1,10]]$  & $630592$\\
   &            & $(8,1)$ & $[[1,8][1,9]]$    & $66064$ \\
\cline{2-5}
   & $L_{10,2}$ & $(2,4)$ & $[[1,2][4,9]]$    & $408216$\\
   &            & $(4,2)$ & $[[1,4][2,9]]$    & $239712$\\
\cline{2-5}
   & $L_{10,3}$ &$(1,1,1,2)$ & $[[2,3][5,8]]$ & $260872$\\
   &            &$(2,1,1,1)$ & $[[2,5][3,8]]$ & $114084$\\
   &            &$(1,2,1,1)$ & $[[3,4][5,7]]$ & $149832$\\
   &            &$(1,1,2,1)$ & $[[3,5][4,7]]$ & $114084$\\
\hline
\end{tabular}
\end{center}
\caption{Values of $\#_{25000}(P)$ for the operators with small
absolute values of the traces.}
\end{table}

In conclusion we formulate the following question. Denote by
$GK(P)$ the probability of the sequence
$P=(a_1,a_2,\ldots,a_{2n-1})$ in the sense of Gauss-Kuzmin:
$$
GK(P)=\frac{1}{\ln(2)}\ln \left(
\frac{(\alpha_1+1)\alpha_2}{\alpha_1 (\alpha_2+1)}\right),
$$
where $\qquad \alpha_1=[a_1{:}a_2;\ldots ;
a_{2n-2};a_{2n-1}],\qquad \alpha_2=[a_1{:}a_2;\ldots ;
a_{2n-2};a_{2n-1}{+}1]$.

\begin{problem}
Let
$$
\begin{array}{ll}
P_1=(a_1,a_2,\ldots a_{2n-1},a_{2n}),&  P_1'=(a_1,a_2,\ldots
a_{2n-1}),\\
P_2=(a_2,a_3,\ldots a_{2n},a_{1}),& P_2'=(a_2,a_3,\ldots a_{2n}).\\
\end{array}
$$
{\it Is the following true}:
$$
\lim\limits_{n\to
\infty}\frac{\#_n(P_1)}{\#_n(P_2)}=\frac{GK(P_1')}{GK(P_2')}\hbox{
?}
$$
\end{problem}


\begin{thebibliography}{99}
\bibitem{Avd1}
M.~O.~Avdeeva, V.~A.~Bykovskii, {\it Solution of Arnold's problem
on Gauss-Kuzmin statistics}, Preprint, Vladivostok,
Dal'nauka,~(2002).
\bibitem{Avd2}
M.~O.~Avdeeva, {\it On statistics of incomplete quotients of
finite continued fractions}, Func. an. and appl., v.~38(2004),
n.~2, pp.~1--11.
\bibitem{ArnProb}
{\it Arnold's problems}, Phasis, M., 2000, Problem 1993-11.
\bibitem{ArnAr} V.~I.~Arnold,
{\it Arithmetics of perfect quadratic forms, symmetry of their
continued fractions and geometry of de Sitter world}, lecture
notes at the conference "Modern Mathematics", Dubna 2002.
\bibitem{Arn2}
V.~I.~Arnold, {\it Continued fractions}, M.: Moscow Center of
Continuous Mathematical Education, (2002).
\bibitem{BSh}
Z.~I.~Borevich, I.~R.~Shafarevich, {\it Number theory}, 3 ed, M.,
(1985).
\bibitem{Hic} D.~R.~Hickerson, {\it Length of period of simple continued
fraction expansion of $\sqrt{d}$}, Pacific J. of Math.,
vol.~46(1973), n.~2, pp.~429-432.
\bibitem{Kar1}
O.~Karpenkov, {\it On tori decompositions associated with
two-dimensional continued fractions of cubic irrationalities},
Func. an. and appl., v.~38(2004), no~2, pp.~28--37.
\bibitem{KarPrep} O.~Karpenkov, {\it Elementary notions of lattice
trigonometry}, preprint 2006,\\
http://arxiv.org/abs/math/0604129
\bibitem{Kat} S.~Katok, {\it Hyperbolic Geometry and Quadratic
forms}, course notes for MATH 497A REU program, summer 2001.
\bibitem{Kle1}
F.~Klein, {\it Ueber einegeometrische Auffassung der gew\"ohnliche
Ketten\-bruchen\-twick\-lung}, Nachr. Ges. Wiss. G\"ottingen
Math-Phys. Kl., 3, (1895), pp.~357--359.
\bibitem{Kle2}
F.~Klein, {\it Sur une repr\'esentation g\'eom\'etrique de
d\'eveloppement en fraction continue ordinaire}, Nouv. Ann. Math.
15(3), (1896), pp.~327--331.
\bibitem{Kor2}
E.~I.~Korkina, {\it Two-dimensional continued fractions. The
simplest examples}, Proceedings of V.~A.~Steklov Math. Ins., v.
209(1995), pp. 143--166.
\bibitem{Lac1}
G.~Lachaud, {\it Voiles et Poly\`edres de Klein}, preprint n
95-22, Laboratoire de Math\'e\-matiques Discr\`etes du C.N.R.S.,
Luminy (1995).
\bibitem{DZg}
J.~Lewis, D.~Zagier, {\it Period functions and the Selberg zeta
function for the modular group}, in The Mathematical Beauty of
Physics, Adv. Series in Math. Physics 24, World Scientific,
Singapore, 1997, pp. 83--97.
\bibitem{Man}
Yu.~I.~Manin, M.~Marcolli, {\it Continued fractions, modular
symbols, and non-com\-mu\-ta\-tive geometry}, 2001,
http://arxiv.org/abs/math/0102006.
\bibitem{Mou2}
J.-O.~Moussafir, {\it Voiles et Poly\'edres de Klein: Geometrie,
Algorithmes et Statistiques},
docteur en sciences th\'ese, Universit\'e Paris IX - Dauphine, (2000),\\
see also at http://www.ceremade.dauphine.fr/$\sim$msfr/
\bibitem{MPav} M.~Pavlovskaya, {\it Continued Fraction Expansions of Matrix
Eigenvectors}, preprint, May 2007.
\end{thebibliography}
\end{document}